\documentclass{amsart}

\newtheorem{theorem}{Theorem}[section]

\theoremstyle{definition}

\theoremstyle{remark}

\DeclareMathOperator{\Char}{char}
\newcommand{\str}{\overrightarrow}
\usepackage{amsfonts}

\numberwithin{equation}{section}
\begin{document}

\title{On $(n+2)$-dimensional $n$-Lie algebras} 
\author{Donald W. Barnes}
\address{1 Little Wonga Rd, Cremorne NSW 2090 Australia}
\email{donwb@iprimus.com.au}
\thanks{This work was done while the author was an Honorary Associate of the
School of Mathematics and Statistics, University of Sydney.}

\begin{abstract} I show that an $(n+2)$-dimensional $n$-Lie algebra over an algebraically closed field must have a
subalgebra of codimension $1$.
\end{abstract}
\maketitle
 R. Bai, X. Wang, H. An and W. Xiao \cite{BWAX} have been working on the classification of the
5-dimensional 3-Lie algebras over an algebraically closed field of characteristic 2.  To complete their classification, they
ask if such an algebra must have a subalgebra of dimension 4.  The following theorem answers that question.

\begin{theorem}  Let $A$ be an $(n+2)$-dimensional $n$-Lie algebra over the algebraically closed field $F$.  Then $A$
has a subalgebra of codimension $1$. \end{theorem}  
\begin{proof}
We denote the derived algebra $[A, A, \dots, A]$ of $A$ by $A^{(1)}$.  If $A^{(1)} < A$, then $A$ has a subalgebra of
codimension $1$ since any subspace containing $A^{(1)}$ is a subalgebra.  Hence we may assume $A^{(1)} = A$ and    
so, that $A$ is not nilpotent.  Let $H$ be a minimal Engel subalgebra.  Then $n-1 \le \dim(H) \le n+1$.  As $F$ is infinite,
$H$ is a Cartan subalgebra by Barnes \cite[Theorem 4.3]{Eng}.  If $\dim(H) = n+1$, the result holds, so we may assume
$\dim(H)\le n$.  This implies that $H$ is abelian and is represented on $A$ by commuting linear transformations.  Since
$F$ is algebraically closed, they have a common eigenvector $u$.  Thus we have $[h_1, h_2, \dots, h_{n-1}, u] =
\alpha(h_1, h_2, \dots, h_{n-1}) u$ for all $h_1, h_2, \dots, h_{n-1} \in H$, where $\alpha$ is a linear map $H^{\wedge
(n-1)} \to F$.  If $\dim(H) = n$, then $\langle H, u \rangle$ is an $(n+1)$-dimensional subalgebra.  

Suppose $\dim(H) = n-1$, $H = \langle a_1, \dots, a_{n-1}  \rangle$.  Let $d$ be the inner derivation $d(a_1, \dots,
a_{n-1})$ of $A$.  For each eigenvalue $\lambda$ of $d$, we have the $\lambda$-component $A_\lambda = \{a \in A
\mid (d-\lambda I)^{n+2} a = 0 \}$ of $A$, where $I$ denotes the identity transformation.  We have $A_0 = H$ and $A$ 
is the direct sum of the components for the eigenvalues of $d$.  Let $\lambda_1, \dots, \lambda_n$ be (not
necessarily distinct)  eigenvalues of $d$.  Then
$[A_{\lambda_1}, \dots, A_{\lambda_n}] \subseteq A_{\lambda_1 + \dots + \lambda_n}$.

Since $H \subset A^{(1)}$, we either have two eigenvalues, say $\alpha, \beta$ with sum $0$ or we have $\alpha + \beta
+
\gamma = 0$.  Suppose first that $\alpha + \beta = 0$.  Suppose $\Char(F) \ne 2$.  We have eigenvectors $u,v$ for
$\alpha, \beta$.  Then $\langle H,u,v \rangle$ is an $(n+1)$-dimensional subalgebra.  Suppose $\Char(F) = 2$.  Then we
can choose $u,v$ such that $[a_1, \dots, a_{n-1},u] = \alpha u$ and $[a_1, \dots, a_{n-1},v] = \alpha v + \theta u$ for 
some $\theta$.  Again we have that $\langle H,u,v \rangle$ is an $(n+1)$-dimensional subalgebra. 

Now suppose $\alpha + \beta + \gamma = 0$.  Suppose $\Char(F) \ne 2$.  Then $\alpha + \beta$ is not an eigenvalue
of $d$, so $[h_1, \dots, h_{n-2}, u,v] = 0$ for all $h_1, \dots, h_{n-2} \in H$.  Thus $\langle H, u,v \rangle$ is an
$(n+1)$-dimensional subalgebra of $A$.  Now suppose $\Char(F) = 2$.  We have the distinct non-zero eigenvalues
$\alpha, \beta, \gamma = \alpha + \beta$ and corresponding eigenvectors $u,v,w$.  If $n = 2,3$, then $H \not\subseteq
A^{(1)}$, so we may suppose $n \ge 4$.

For some re-ordering of the basis of $H$, we have $[a_1, \dots, a_{n-3}, u, v, w] \ne 0$.  Denote the string $a_1, \dots,
a_{n-3}$ by $\str{a}$.  We apply the Jacobi identity to  the product $P = [\str{a}, a_{n-2}, u, [\str{a}, a_{n-1}, v,w]]$.  Since
$[\str{a}, a_{n-1}, v,w] \in \langle u \rangle$, $P = 0$.  But
\begin{equation*}\begin{split}
P &= [\str{a},  [\str{a}, a_{n-2}, u, a_{n-1}], v,w] + [\str{a}, a_{n-1}, [\str{a}, a_{n-2},  u, v], w]\\
& \qquad\qquad \qquad\qquad+ [\str{a}, a_{n-1}, v, [\str{a},
a_{n-2}, u, w]]\\
&= [\str{a}, \alpha u, v,w] + 0 + 0\\
\end{split}\end{equation*}
since $ [\str{a}, a_{n-2},  u, v] \in \langle w \rangle$ and $ [\str{a}, a_{n-2}, u, w] \in \langle v \rangle$.  Therefore
$\alpha = 0$ contrary to the  definition of $\alpha$.  Thus this case cannot arise.
\end{proof}
\bibliographystyle{amsplain}

\end{document}